\documentclass[12pt]{article}
\usepackage{amscd,amsmath,amsfonts,amssymb,amsthm,amsxtra,indentfirst,mathrsfs}

\usepackage[all]{xy}

\newtheorem{prop}{Proposition}[section]
\newtheorem{thm}[prop]{Theorem}
\newtheorem{lem}[prop]{Lemma}

\theoremstyle{definition}

\newtheorem{defn}[prop]{Definition}
\newtheorem{rem}[prop]{Remark}
\newtheorem{exa}[prop]{Example}

\newtheorem{conv}[prop]{Convention}

\newcommand{\al}{\alpha}
\newcommand{\be}{\beta}

\newcommand{\de}{\delta}
\newcommand{\ep}{\epsilon}

\newcommand{\la}{\lambda}
\newcommand{\om}{\omega}

\newcommand{\Om}{\Omega}

\newcommand{\cA}{\mathcal A}
\newcommand{\cB}{\mathcal B}
\newcommand{\cC}{\mathcal C}
\newcommand{\cD}{\mathcal D}
\newcommand{\cF}{\mathcal F}
\newcommand{\cG}{\mathcal G}
\newcommand{\cH}{\mathcal H}

\newcommand{\cT}{\mathcal T}

\newcommand{\cR}{\mathcal R}

\newcommand{\bC}{\mathbb C}

\newcommand{\bR}{\mathbb R}
\newcommand{\bZ}{\mathbb Z}

\newcommand{\fS}{\mathfrak S}

\newcommand{\tP}{\tilde{P}}
\newcommand{\ttP}{\tilde{\tilde{P}}}
\newcommand{\tB}{\tilde{B}}
\newcommand{\tb}{\tilde{b}}

\newcommand{\qq}{\qquad\quad}
\newcommand{\q}{\quad}

\author{Nir Avni}

\title{Entropy Theory for Cross Sections}

\begin{document}

\maketitle

\abstract{We define the notion of entropy for a cross-section of an action of continuous amenable group, and relate it to the entropy of the ambient action. As a result, we are able to answer a question of J.P. Thouvenot about completely positive entropy actions.}

\section{Introduction}
\subsection{Cross Sections of Flows}
One of the oldest techniques for studying continuous flows is to restrict them to cross-sections. For a smooth flow---that is, an action of the group $\bR$ by diffeomorphisms on a smooth manifold $M$---Poincare defined a cross-section as a co-dimension 1 sub-manifold, $N$, that is transverse to the direction of the flow. The transversality condition implies that each orbit intersects $N$ in a discrete set. In other words, for each point $x$ of $M$, the set of visit times, $\{t\in\bR | T_t(x)\in N\}$ is discrete. The co-dimension 1 condition implies that $N$ intersect many orbits. The Poincare return map assigns to a point $x\in N$ the first point in the positive trajectory of $x$ that lies in $N$.

Extending this construction, Ambrose proved (see \cite{Am}) that for any Borel flow on a Borel probability space $(X,\mu)$, there is a Borel subset $Y\subset X$ such that, for almost all points $x\in X$, the set $\{t\in\bR_{>0}|T_t(x)\in Y\}$ is non-empty and discrete. If the flow preserves the measure $\mu$, then there is a measure $\nu$ on $Y$ such that, near $Y$, the measure $\mu$ is a product of $\nu$ and the Lebesgue measure (in the direction of the flow). It follows that for $\nu$-almost all $y\in Y$, the set $\{t\in\bR_{>0}|T_t(y)\in Y\}$ is also nonempty and discrete. Therefore, this set has a minimum, which we call $\al (y)$. Define a transformation $S:N\to N$ by $S(y)=T_{\al (y)}(y)$. One can show that $S$ preserves the measure $\nu$. 

In this way, a flow $(X,T_t)$ induces a measure preserving transformation $(Y,S)$. Conversely, Ambrose shows that the flow can be obtained from the triple $(Y,S,\al)$ by the construction of the `flow under a function'. There are many connections between the dynamics of the flow $T_t$ and the dynamics of the transformation $S$; one of the basic ones is the Abramov formula (see \cite{CFS}) connecting the entropies of the two actions:
\begin{equation} \label{eq:Abramov.one.dim}
h(S)=h(T_1)\int \limits _{Y}\al (y) d\nu .
\end{equation}

The method of cross-sections is less-understood for higher dimensional groups. In \cite{FHM}, the authors show that every Borel action of a locally-compact and second-countable group admits a cross-section (see Definition \ref{defn:cross.section} below). However, even for the group $\bR^2$, the cross-section does not come equipped with an action of a discrete subgroup of the continuous group. In this case, it is not clear what is the higher dimensional analogue of Abramov formula.

\subsection{Main Results}
The purpose of this article is to develop the technique of cross-sections for probability-preserving free actions of a class of amenable groups. This class contains, for example, all nilpotent Lie groups. In particular, we state and prove a generalization of Abramov formula for these groups. As a corollary, we deduce a generalization of a theorem of Rohlin and Sinai on Kolmogorov systems. Before stating the theorem, we give some definitions. Let $(X,\cB ,\mu )$ be a probability space. For a partition $P$ of $X$, whose parts have measures $p_1,\ldots ,p_n$, let $H(P)=\sum p_i \log p_i$. A transformation $T:X\to X$ is called {\em uniformly mixing} if, for every partition $P$ and $\ep >0$, there is an integer $N$ such that for every $k$ and every finite sequence of integers $i_1 < i_2 < \dots < i_k$ satisfying $i_{j+1}-i_j \geq N$ for all $j$, we have
\[
\left | \frac{1}{k}H(T^{i_1}P \vee T^{i_2}P \vee \dots \vee T^{i_k}P) - h(T,P) \right | <\ep .
\]
Here $h(T,P)$ is the (measure theoretic) entropy of $T$ with respect to the partition $P$.

Recall that the {\em spectrum} of $(X,\cB ,m,T)$ is defined as the spectrum of the unitary operator $U_T :L_2(X,\cB ,m)\to L_2(X,\cB ,m)$ given by $(U_T f)(x)=f(Tx)$.

In \cite{RS}, the following theorem is proved:

\begin{thm} \label{Rohlin:Sinai} Let $(X,\cB ,m)$ be a probability space and let $T:X\to X$ be a probability-preserving transformation. Assume that, for any nontrivial partition $P$ of $X$, we have $h(T,P)>0$. Then
\begin{enumerate}
\item The transformation $T$ is uniformly mixing.
\item The spectrum of $(X,\cB ,m,T)$ is Lebesgue with countable multiplicity.
\end{enumerate}
\end{thm}

An entropy theory for actions of amenable groups is developed in \cite{OW}. The definitions there do not apply to all amenable groups, but only to a class of amenable groups called {\em groups with good entropy theory}. Let $G$ be an amenable group with good entropy theory. Given a probability-preserving action of $G$ on a probability space $X$ and a partition $P$ of $X$, the entropy of the action with respect to the partition will be denoted by $h(G,P)$.

The entropy of an action relative to a sub-sigma-algebra is defined only for a smaller class of amenable groups. These are called in \cite{OW} {\em groups with zero self entropy}. We note here that all nilpotent Lie groups have zero self entropy, and that all groups with zero self entropy (and, in fact, all groups with good entropy theory) are unimodular. For more information on groups with zero self entropy, see Section \ref{sec:ent.theory}. Starting from Section \ref{sec:ent.theory}, we shall assume that the groups we talk about have zero self entropy.

For technical reasons, we deal with Borel actions only. This means that the probability space $X$ is a Polish space endowed with the Borel sigma-algebra, that the probability measure is regular, and that the action map $G\times X \to X$ is Borel measurable. In fact, as remarked in \cite{FHM}, since we are interested only in measure theoretic properties, we can change $X$ and the action to an isomorphic (in the category of measure-preserving actions) Borel action, which is also continuous.

An action of $G$ on $X$ is called {\em completely positive entropy} (or CPE for short) if, for any non-trivial partition $P$ of $X$, the entropy $h(G,P)$ is strictly positive. It is called {\em free} if, for every $g\in G$, which is different from $1$, the set of fixed points of $g$ is negligible.

\begin{defn} Let $G$ be a group and let $K\subset G$ be a subset. A finite subset $F\subset G$ is called $K$-separated if, for every two non-equal elements $g,h$ of $F$, the element $gh^{-1}$ does not belong to $K$.
\end{defn}

This is the generalization of Theorem \ref{Rohlin:Sinai}:

\begin{thm} \label{thm:cpe} Let $G$ be an amenable group with zero self entropy. Suppose $G$ acts freely on a probability space $X$, and suppose that the action is CPE. Then
\begin{enumerate}
\item For every partition $P$ and any $\ep>0$, there is a compact set $K\subset G$ such that for any finite set $F\subset G$ that is $K$-separated,
\[
\left | \frac{1}{|F|}H\left ( \bigvee \limits _{g\in F} gP \right ) -H(P)\right | < \ep .
\]
\item As a $G$-module, the space $L_2(X)$ decomposes as a direct sum of infinitely many copies of the regular representation of $G$.
\end{enumerate}
\end{thm}

This theorem was conjectured by Thouvenot. It was proved for discrete amenable groups in \cite{RW} and \cite{DG}; see also \cite{D} and \cite{DP}.
\subsection{Organization}
In the rest of this introduction, we shall outline, without proofs, the main points of this paper.

Let $G$ be a locally compact, second countable, and unimodular group. Let $X$ be a Polish space, let $\cB$ be the Borel sigma-algebra on $X$, and let $m$ be a Borel probability measure on $X$. For a free, Borel, and measure-preserving action of $G$ on $(X,\cB ,m)$, we define a {\em cross-section} to be a Borel set $S\subset X$ that intersects almost every orbit in a discrete set. On $S$, we have the Borel sigma-algebra $\cB _S$. We shall show that there is a canonically defined Borel measure $\mu$ on $S$, such that, locally near $S$, the measure $m$ is the product of $\mu$ and the Haar measure of $G$ (in the direction of the action). 

In contrast to the one-dimensional case, $S$ does not come equipped with a canonical action of a discrete group. There are, however, additional structures on $S$. First, there is an equivalence relation: $x\sim y$ if and only if $x$ and $y$ (both are elements of $S$) lie in the same $G$ orbit. Denote by $\cR \subset S\times S$ the set of equivalent pairs. We also get a function $\al : \cR \to G$ by defining $\al (x,y)=g$ if $gy=x$ (this $g$ is unique because the action is assumed to be free). Note that we used the same letter, $\al$, as before, since it is an analogue of the function we introduced in the one-dimensional case.

We concentrate on quintuples of the form $(S,\cB ,\mu ,\cR ,\al )$. Note that this information contains $G$ implicitly. Given $\fS= (S,\cB ,\mu ,\cR ,\al )$, it is possible to formulate an analogue of the mean ergodic theorem. If $\fS$ comes from a cross-section of an action of $G$, then this analogue of the mean ergodic theorem holds. Surprisingly, the converse is also true: the mean ergodic theorem for $\fS$ implies that $\fS$ is isomorphic to a cross-section of a probability-preserving action. These constructions are described in Section \ref{sec:erg.theory}. In the rest of the article, we only deal with such quintuples, which we call cross-sections (slightly abusing notations). In Section \ref{sec:erg.theory}, we also describe our main technical tools---a tiling lemma that generalizes \cite[Proposition 7]{OW}, and an ergodic theorem for (generalizations of) Rohlin towers.

Section \ref{sec:ent.theory} is devoted to the entropy theory of cross-sections. Building on \cite{OW}, we define entropy for cross-sections of actions of groups with zero self entropy, and prove an analogue of Abramov's theorem. Interestingly, the useful notion here is relative entropy (with respect to a sub-sigma-algebra) rather than the absolute entropy.

In the same section, we prove the following transfer theorem (see also \cite{RW}): Let $G,H$ be two amenable groups with zero self entropy, let $\fS =(S,\cB ,\mu ,\cR ,\al )$ be a cross-section for an action of $G$, and let $\fS '=(S', \cB ',\mu ',\cR ',\be )$ be a cross-section for an action of $H$. Suppose that $\phi :S\to S'$ is a measure-preserving and equivalence-preserving map (i.e. $x\cR y \iff \phi (x)\cR ' \phi (y)$), and assume that $\be \circ \phi$ is measurable with respect to a sub-sigma-algebra $\cG \subset \cB$. Then, the relative entropies of $\fS$ and $\fS '$ with respect to $\cG$ and $\phi _* \cG$ are the same. 

This transfer theorem is useful to the ergodic theorem of groups. By a theorem of \cite{CFW}, for any $\fS$ and $\fS '$ as above, there is an isomorphism $\phi$ as above (in general, however, we cannot say much about $\cG$). This allows us to transfer questions (and answers) from actions of one group to actions of another, using cross-sections as intermediaries.

In Section 4, we prove Theorem \ref{thm:cpe}. Section 5 features some concluding remarks.
\subsection{Acknowledgment}
This work contains results of a research done under the supervision of Prof. Benjamin Weiss at the Hebrew University. I thank him for being such a good source of interesting questions, for so many illuminating discussions, and for his help in the process of writing this paper. I also thank Prof. Dan Rudolph, Prof. Valentyn Golodets, and the referees for helpful remarks.

\section{Ergodic Theory for Cross Sections} \label{sec:erg.theory}

\subsection{Cross Sections}

\begin{conv} Unless stated otherwise, all groups will be locally compact, second countable, unimodular, and amenable. Starting from the next section, we also assume that they have zero self entropy. All spaces are assumed to be Borel, all sigma-algebras are assumed to be contained in the Borel sigma-algebra, and all actions are assumed to be continuous, probability-preserving, and free.
\end{conv}

Suppose $G$ is a group as above with Haar measure $\la$. Recall that $G$ is called {\em amenable} if, for every $\ep >0$ and a compact $K\subset G$, there is a compact $F\subset G$ such that $\la (KFK)<(1+\ep )\la (F)$. Such an $F$ is called {\em $(K,\ep )$-invariant}. We say that a compact set $F$ is {\em sufficiently invariant} if $F$ is $(K,\ep )$-invariant for some $K$ and $\ep$ (which should be prescribed). A sequence $(F_n)$ of compact subsets of $G$ is called a {\em F\o lner sequence} if, for every $K$ and $\ep$, there is an $N$ such that $F_n$ is $(K,\ep )$-invariant for $n>N$. We refer the reader to \cite{OW} and \cite{W} for the ergodic theory of actions of amenable groups.

If $G$ acts on a probability space $(X,\cB ,m)$, we denote the application of the group element $g$ to the point $x$ by $g\cdot x$.

\begin{defn} \label{defn:cross.section} Suppose that the group $G$ acts freely on a probability space $(X,\cB,m)$. A Borel subset $S\subset X$ is called a {\em cross-section} for the action if, for almost every point $x\in X$, the set of return times of $x$---which is defined as $\{ g\in G | gx \in S\}$---is discrete and non-empty.
\end{defn}
By a theorem of \cite{FHM}, we know that cross-sections of free Borel actions of local compact, second countable groups always exist.

We will usually work with a slightly stronger condition. For a neighborhood $U \subset G$ of the identity, we say that the cross-section is $U$-discrete if, for almost any $x\in X$, the set of return times of $x$ is $U$-separated. In an ergodic system, for any compact neighborhood $U$ of the identity in $G$, every cross-section contains a subset which is a cross-section and is $U$-discrete.

Clearly, if $S$ is a cross-section, then the restriction of the Borel sigma-algebra to $S$ is the Borel sigma-algebra on $S$. We denote it by $\cB$.

\begin{defn} For a $U$-discrete cross-section $S$, the induced measure on $S$ is
\[
\mu (A)=\lim \limits _{V \searrow \{ 1\} }\frac{m(V\cdot A)}{\la (V)}
\]
where the $V$'s are neighborhoods of $1$ in $G$. The limit over $V \searrow \{ 1\}$ means that for every open neighborhood $W$ of $1$, all but finitely many of the $V$'s are contained in $W$.
\end{defn}
If the cross-section is $U$-discrete and $V^2 \subset U$, then we have a one-to-one map $a:V\times S \to X$. Pulling back the measure $m$ on $V\cdot S$ gives a measure $a^* m$ on $V\times S$. We decompose $a^* m$ with respect to $S$: $a^* m =\int _S \eta _x d\mu (x)$. Since the action is measure-preserving, almost every $\eta _x$ must be a multiple of the Haar measure $\la$ of $G$. After multiplying $\mu$ by some function we can assume that $\eta _x =\la$ for almost every $x$. It is now clear that the limit above exists (it does not even depend on $V$ if $V$ is small enough). Also, $\mu (S)=m(V\cdot S)/\la (V)$, so $\mu$ is a finite measure.

\begin{conv} \label{intensity1} By rescaling the Haar measure, we can (and will) assume $\mu (S)=1$.
\end{conv}

The last two ingredients we need are the equivalence relation and the cocycle induced by the action. The first is the subset $\cR \subset S\times S$ that consists of pairs $(x,y)$ in the same $G$ orbit. If $(x,y)\in \cR$ we also write $x\cR y$. Clearly, $\cR$ is a Borel set and is an equivalence relation. If $(x,y)\in \cR$, then by definition there is a $g\in G$ such that $g\cdot y=x$. This $g$ is unique because we assumed the action was free. We set $\al (x,y)=g$. The function $\al$ is measurable and satisfies the equation 
\[
\al (x,y) \al (y,z) =\al (x,z).
\]

Functions from an equivalence relation to a group satisfying this equation are called cocycles\footnote{The reader may be more familiar with the definition of cocycle of a transformation $T:X\to X$ as being a function $f:\bZ \times X \to \bR ^{\times}$ that satisfies $f(n+m,x)=f(n,x)\cdot f(m,T^nx)$. This is a special case of our definition, where $G$ is replaced by $\bR ^{\times}$, the equivalence relation consists of pairs of the form $(x,T^n x)$ and $\al (T^n x,x)=f(n,x)$.}. We call $\al$ the {\em cocycle induced by the action}. The cocycle $\al$ that we have constructed satisfies the condition that for $\mu$-almost all $x\in S$, the function $\al(x,-)$ (from the equivalence class of $x$ to $G$) is one-to-one. We call cocycles that satisfy this condition {\em free}. Moreover, if the cross-section $S$ is $U$-discrete, then for every $x,y\in S$ such that $x\cR y$ and $x\neq y$, the element $\al(x,y)$ is not in $U$. A cocycle that satisfies this property is called $U$-discrete.

The following lemma is well-known.

\begin{lem} Let $G$ be a unimodular group, let $U\subset G$ be an open neighborhood of the identity, and let $X$ be a probability space. Suppose that $S$ is a $U$-discrete cross-section of a probability-preserving action of $G$ on $X$. For every Borel transformation $T:S\to S$ such that $x\cR(Tx)$ for almost any $x$, and for any $A\subset S$ such that the restriction of $T$ to $A$ is one-to-one, we have $\mu (TA)=\mu (A)$ (in this situation we say that $\cR$ preserves the measure $\mu$).
\end{lem}

\begin{proof} Since $T|_A$ is one-to-one, it is enough to prove the lemma for every element in a (countable) decomposition $A=\cup A_i$ of $A$. Let $V$ be a neighborhood of $1$ such that $V^4 \subset U$. The inverse of the action map is a measure-preserving isomorphism $\Phi :V^2\cdot S \to V^2\times S$ (where on the left we take the measure $m$ and on the right we take $\la \times \mu$). Consider the function $x\mapsto \al (x,Tx)$. There is a decomposition of $A$ into sets $A_i$ such that for every $x,y$ in the same $A_i$,
\begin{equation} \label{eq:al.al.-1}
\al (x,Tx) \al (y,Ty) ^{-1} \in V.
\end{equation}
By decomposing $A$ to the $A_i$, we can assume that (\ref{eq:al.al.-1}) holds for every two points in $A$. Fix some $x_0$ in $A$, and let $g=\al (x_0 ,Tx_0)$. Let $W$ be a neighborhood of $1$ in $G$ such that $W,g^{-1}Wg\subset V$. Since the action of $G$ preserves $m$ and $G$ is unimodular,
\[
\mu (A)=\frac{m(g^{-1}Wg\cdot A)}{\la (g^{-1}Wg)} = \frac{m(Wg\cdot A)}{\la (W)} = \frac{(\la \times \mu )(\Phi (Wg\cdot A))}{\la (W)}.
\]
For every $x\in A$, $\Phi (Wg\cdot x)$ is of the form $Wh \times \{ Tx\}$ where $h\in V$ is defined by $hTx = gx$. Hence $(\la \times \mu )(\Phi (Wg\cdot A))=\la (W) \mu (TA)$, so $\mu(A)=\mu(TA)$.
\end{proof}
\subsection{The Ergodic Theorem}
From now on we shall work with quintuples $(S,\cB ,\mu ,\cR ,\al )$ where $(S,\cB )$ is a Borel space, $\mu$ is a probability measure on $S$, $\cR$ is a measure-preserving equivalence relation on $S$, and $\al$ is a free cocycle with values in an amenable group $G$. We stress again that $G$ is implicitly given in this data. We wish to develop ergodic theory of such quintuples. The first step is a mean ergodic theorem:

\begin{thm} \label{mean:erg:thm} Let $(F_n)$ be a F\o lner sequence in $G$. Suppose that $(S,\cB ,\mu ,\cR ,\al )$ is a quintuple corresponding to a cross-section of a probability-preserving free action of $G$, and suppose that $S$ is $U$-discrete for some open neighborhood $U$ of the identity in $G$. Then for any $h\in L_{\infty}(S,\mu )$, the sequence of functions
\[
h_n (x)=\frac{1}{|\{ (x,y)\in \cR\wedge \al (x,y)\in F_n \}|}\sum \limits _{(x,y)\in \cR\wedge \al (x,y)\in F_n} h(y)
\]
(and $h_n(x)=0$ if the set $\{y | \al(x,y)\in F_n\}$ is empty) converges in measure to the constant function $\int _S h d\mu$.
\end{thm}

\begin{proof} Let $M$ be such that $h(x)\leq M$ for almost every $x$. Let $V\subset G$ be a neighborhood of the identity such that $V^2\subset U$.  Define a function $H$ on $X$ by 
\[
H(x)=\left \{ \begin{matrix} h(s) & \qq & {(\exists g\in V , s\in S)\quad x=gs } \\ 0 & \qq & \textrm{Otherwise} \end{matrix} \right . .
\]
First note that for any compact subset $F\subset G$ and for almost any $s\in S$,
\[ 
\left | \sum _{t: \al (s,t)\in F} h(t)-\frac{1}{\la (V)}\int _{F}H(fs)df \right | \leq M\frac{\la (VF\setminus F)}{\la (V)} .
\] 
If $n$ is big enough then $F_n$ is $(V,\ep )$-invariant, so $\la(VF_n\setminus F_n)<\ep\la(F_n)$. By the mean ergodic theorem for the action of $G$ on $X$, we have that for any $\ep >0$, if $n$ is big enough then for all $x\in X$ outside a set of measure less than $\ep \la (V)$,
\[ 
\left | \frac{1}{\la (F_n)} \int _{F_n} H(fx)df - \la (V)\int _S h(s)ds \right | = \left | \frac{1}{\la (F_n)}\int _{F_n} H(fx)df - \int _X H(y)dy \right | <\ep .
\]
By Fubini's theorem, there is a $g\in V$ such that the above inequality holds for all $x\in gS$ except a subset of $g_* \mu$-measure less than $\ep$. Since we can assume that $F_n$ is $(V, \ep )$-invariant, we have that for {\em all} $x\in X$,
\[
\frac{1}{\la (F_n)}\left | \int _{F_n} H(fx)-H(fgx)df\right | \leq \frac{2}{\la (F_n)}\int _{F\triangle Fg} |H(fx)|df \leq 2\ep M,
\]
so for a subset of $S$ of $\mu$-measure greater than $1-\ep$ we have
\[
\left | \frac{1}{\la (F_n)} \sum _{t: \al (s,t)\in F_n} h(t) - \int _S h(t)dt \right | \leq \frac{1}{\la (F_n)\la (V)} \left |  \int _{F_n} H(fs)df - \int _{F_n}H(fgs)df\right | + 
\]
\[
+ \frac{(M+1)\ep}{\la (V)} \leq \frac{(2M+2)\ep}{\la (V)}.
\]
Taking $\ep$ small enough, we get that for all bounded functions $h$,
\[
\frac{1}{\la (F_n)}\sum _{t:\al (s,t)\in F_n}h(t) \to \int _S h(t)dt
\]
in probability. Applying this to the function $h=1$, we get that outside a set of small measure,
\[
1-\ep \leq \frac{|\{ t: \al (s,t)\in F_n \} |}{\la (F_n)} < 1+\ep.
\]
Dividing the last two inequalities gives the result.
\end{proof}

\begin{rem} In a similar way, one can prove a pointwise ergodic theorem (along a tempered F\o lner sequence) by using \cite{L} instead of the mean ergodic theorem. A direct argument can also be given using Proposition \ref{prop:tiling}.
\end{rem}

Conversely, the mean ergodic theorem implies that the quintuple comes from a $G$ action on a probability space:

\begin{thm} \label{inverse:mean:erg:thm} Suppose $(X,\cB ,\mu ,\cR ,\al )$ is a quintuple where $(X,\cB,\mu)$ is a probability space, $\cR$ is a Borel equivalence relation with countable equivalence classes that preserves $\mu$, and $\al:\cR\to G$ is a free cocycle whose range is some amenable group $G$. Assume that for almost all $x\in X$, the set $\{\al(x,y)|(x,y)\in\cR\}$ is discrete, and that for any $f\in L_{\infty}(X,\mu )$, the sequence of functions $f_n$ defined in the last theorem converges in probability to $\int fd\mu$. Then there is a probability-preserving action of $G$ such that the quintuple is induced from this action.
\end{thm}

\begin{proof} We take the Mackey range of the cocycle. That is, we look at $X\times G$ and divide by the equivalence relation $(x,g)\sim (y,\al (x,y)g)$ for all $(x,y)\in \cR$ and $g\in G$. According to \cite{FHM} the quotient $Y=X\times G /\sim$ is a Borel space. 

The group $G$ acts on $Y$ by right multiplication. Since we assume that $\al$ is free, the map $x\mapsto [(x,1)]$ is an isomorphism. Let $S\subset Y$ be the set $\{[(x,1)]\}$. It follows from the assumptions that almost every $G$-orbit in $Y$ intersects $S$ in a discrete set. Since every $G$-orbit intersects $S$ at least once, we get that $S$ is a cross-section. The relation $\sim$ on $X\times G$, as well as the right action of $G$ preserve the measure $\mu \times \la$, and thus the quotient inherits a $G$-invariant measure. It remains to show that the mean ergodic theorem implies that this measure is finite.

Let $\pi :X\times G \to Y$ be the quotient map. We show that if $F\subset X\times G$ is Borel, its boundary has measure zero, and such that $\pi | _F$ is one-to-one, then $\mu \times \la (F)\leq 1$ (recall that the measure on the quotient can be computed from a Borel fundamental domain $Z$ with negligible boundary by $m(A)=\mu \times \la (\pi ^{-1}(A))\cap Z$ for $A\subset X\times G / \sim$). Suppose $\mu \times \la (F) >1+2\ep$. After removing from $F$ a subset of measure less than $\ep$, we can assume that there is a partition $X=X_1 \sqcup \dots \sqcup X_n$ and compact subsets $K_1 ,\dots ,K_n$ such that $F=\cup X_i \times K_i$. Let $\eta$ be a positive small number (to be determined later), and let $L$ be a $(\cup K_i ,\eta)$-invariant set such that all but $\eta$ of the points in $X$ satisfy the mean ergodic theorem with respect to all the characteristic functions $1_{X_i}$.

Let $s\in X$ be such a point. If $t\in Ls\cap X_i$ (we shall write this as $i=i(t)$), then $K_i\al (s,t)\subset (\cup K_j) L$, and the sets $K_{i(t)}\al (s,t)$, for the different $t$'s, are disjoint. Therefore,
\[
(1+\eta )\la (L) \geq \sum _{t\in Ls\cap X} \la (K_{i(t)}\al (s,t)) = \sum _i \la (K_i) |\{ t| t\in Ls\cap X_i\}| \geq
\]

\[
\geq \sum _i \la (K_i) \mu (X_i) \la (L) (1-\eta ) = (1-\eta )\la (L) (\mu \times \la )(F),
\]
which is a contradiction if $\eta$ is small enough.
\end{proof}

\begin{conv} From now on, by a cross-section we will mean a quintuple $(S,\cB ,\mu ,\cR ,\al )$---where $(S,\cB,\mu)$ is a Borel probability space, $\cR$ is a Borel equivalence relation with countable equivalence classes that preserves $\mu$, and $\al:\cR\to G$ is a free cocycle which is $U$-discrete---that satisfies the mean ergodic theorem.
\end{conv}

\subsection{Tiling Lemma}
\begin{defn} A sequence of finite sets $A_i \subset G$ is called {\em $\ep$-disjoint} if, for every $i$,

\[
\left|A_i \cap \bigcup \limits _{j<i} A_j\right|<\ep |A_i| .
\]
\end{defn}

\begin{lem} \label{lem:tiling} Let $G$ be an amenable group, let $c$ be an integer, and let $0<\delta <0.1$. Let $U$ and $V$ be neighborhoods of the identity in $G$ such that $V^2\subset U$, and let $F\subset G$ be a compact set such that there are at most $c\lambda (F)$ disjoint right translates of $V$ whose centers are in $F$. Assume also that $\lambda (F)>10$. Let $A\subset G$ be a finite set which is $U$-separated, and let $B\subset A$ be such that for every $b\in B$,
\[
|Fb\cap A|>\frac{1}{2}\lambda(F).
\]
Then there is a subset $\tB =\{\tb _1,\tb _2,\ldots,\tb _k\}\subset B$ such that the sequence of sets $F\tb_i\cap A$, $i=1,\ldots,k$ is $\delta$-disjoint, and such that
\[
|F\tB \cap A| > \frac{\delta}{4c}|B| .
\]
\end{lem}

\begin{proof} Let $B=\{ b_1,b_2,...,b_{|B|}\}$. Let $\tb _1=b_1$, and for $i>1$ let $\tb _i$ to be the first element of $B$ such that
\[
|F\tb _i\cap \bigcup \limits _{j<i} F\tb _j \cap A|<\delta |F\tb _i \cap A|,
\]
if such a $\tb$ exists. Otherwise stop, and let $k=i-1$. Let $\tB =\{ \tb _1,...,\tb _k \}$. One of the following must hold:\\
1. $|\tB|>|B|/2$. In that case, every $\tb _i\in \tB$ contributes at least $\frac{1-\delta}{2}\lambda (F)$ new points to the union $F\tB\cap A$, so
\[
\left|\bigcup \limits _i F\tb _i \cap A\right| > \frac{1-\de}{4}\lambda (F)|B|>|B|.
\]
2. $|\tB|<|B|/2$. Let $E=B\setminus\tB$. For every point $b\in E$, more than a $\de$ portion of the set $Fb\cap A$ is contained in $F\tB\cap A$. This implies that for every $b\in E$,
\[
\langle 1_{Fb\cap A},1_{F\tB \cap A}\rangle \geq \delta |Fb\cap A|
\]
where the inner product of the functions is taken with respect to counting measure. Since by our assumptions, the function $\sum \limits _{b\in B}1_{Fb\cap A}$ is smaller than $c\lambda (F)$,
\[
c\lambda (F) |F\tB \cap A| \geq \langle \sum \limits _{b\in E} 1_{Fb\cap A},1_{F\tB \cap A}\rangle \geq \delta \sum \limits _{b\in E} |Fb\cap A| \geq \frac{1}{4} \delta \lambda (F)|B|
\]
 Thus, in both cases the conclusion of the lemma holds.
\end{proof}

\begin{prop} \label{prop:tiling} Let $G$ be an amenable group, and let $c$ be an integer. For every $0<\delta <0.1$ there are $\epsilon >0$ and an integer $N$, such that if $U,V\subset G$ are neighborhoods of the identity, $F_1,...,F_N$ are compact subsets of $G$, and $B\subset A \subset G$ are finite sets such that 
\begin{enumerate}
\item $V^2\subset U$.
\item For every $i$ there are at most $c\lambda (F_i)$ disjoint translates of $V$ whose centers are in $F_i$.
\item $\lambda (F_i)>10$.
\item $A$ is $U$-separated.
\item $|B|>(1-\epsilon)|A|$.
\item For every $i$ and every $b\in B$, we have $\left|\left(\bigcup\limits _{j< i}F_j\right)^{-1} F_ib\cap A\right| < (1+\epsilon)|F_ib\cap A|$.
\end{enumerate}

Then there are subsets $\tB _i =\{ \tb _{i1},...,\tb _{ik_i} \} \subset B$, $i=1\ldots N$ such that 
\begin{enumerate}
\item For any $i$, the sequence of sets $F_i\tb _{i1}\cap A,...,F_i\tb _{ik_i}\cap A$ is $\delta$-disjoint.
\item For any $i_1 \neq i_2$, the sets $F_{i_1}\tB _{i_1}\cap A$ and $F_{i_2}\tB _{i_2}\cap A$ are disjoint.
\item $\left|\bigcup \limits _i F_i\tB _i\cap A\right| >(1-2\delta )|A|$.
\end{enumerate}
\end{prop}

\begin{defn} We say that a sequence $(A,B,F_i)$ satisfies the conditions of the tiling lemma with parameter $\epsilon$ if the above conditions hold.
\end{defn}

\begin{proof} Let $N=\frac{log(\delta )}{log(1-\delta /8c)}$, and let $\epsilon>0$ be such that $(1-\delta )(1-2\epsilon )>(1-2\delta)(1+\delta+2\epsilon )$---note that when $\ep=0$, the left hand side is equal to $1-\de$ and the right hand side is equal to $(1-2\de)(1+\de)=1-\de-2\de^2$, which is strictly smaller than the left hand side.

We define three sequences of sets, $A_i,B_i$, and $\tB_i$ for $0<i\leq N$. Let $A_N=A$, $B_N=B$. Apply Lemma \ref{lem:tiling} to $(A_N,B_N,F_N)$ to get a set $\tB_N$ such that the collection of sets $\{F_Nb\cap A_N=F_Nb\cap A | b\in\tB_N\}$ is $\de$-disjoint, and the size of its union is bigger than $\frac{\de}{4c}|B_N|$. Assume that $A_i,B_i$ have been defined for $k\leq i\leq N$. If $|A_k|<2\de|A|$, we define $A_{k-1}=A_k,B_{k-1}=B_k,\tB_{k-1}=\emptyset$. Otherwise, let 
\[
A_{k-1}=A_{k}\setminus F_k\tB _k,
\]
and
\[
B_{k-1}=B_k\setminus \left ( \bigcup \limits _{j<k}F_j \right ) ^{-1}F_k\tB _k,
\]
and apply Lemma \ref{lem:tiling} to $(A_{k-1},B_{k-1},F_{k-1})$ to get a set $\tB_{k-1}\subset B_{k-1}$, such that the collection of sets $\{F_{k-1}b\cap A_{k-1}=F_{k-1}b\cap A | b\in\tB_{k-1}\}$ is $\de$-disjoint, and the size of its union is bigger than $\frac{\de}{4c}|B_{k-1}|$. The first two conclusions of the proposition hold. For the third, we use

\begin{lem} For every $k$, if $|A_k|>2\de|A|$, then $|B_k|>|A_k|/2$.
\end{lem}

\begin{proof} By decreasing induction on $k$. The case $k=N$ follows from assumption $5$. Assuming the lemma for $k+1,\ldots,N$, let
\[
\Sigma=\sum_{j=k+1}^N\sum_{b\in\tB_j}|F_jb\cap A|
\]
and
\[
\Sigma'=\sum_{j=k+1}^N\sum_{b\in\tB_j}\left|\left(\cup_{i<j}F_i\right)^{-1}F_jb\cap A\right|.
\]
By $\de$-disjointness and the assumption of the lemma,
\begin{equation} \label{ineq:tiling}
2\de|A|<|A_k|\leq|A|-(1-\de)\Sigma \Longrightarrow \Sigma < \frac{1-2\de}{1-\de}|A|.
\end{equation}
By assumption 6,
\[
|B_k|\geq|B|-\Sigma'\geq|B|-(1+\ep)\Sigma.
\]
Hence, by assumption 5,
\[
\frac{|B_k|}{|A_k|}\geq\frac{|B|-(1+\ep)\Sigma}{|A|-(1-\de)\Sigma}\geq \frac{(1-\ep)|A|-(1+\ep)\Sigma}{|A|-(1-\de)\Sigma}.
\]
We claim that the last expression is larger than $1/2$. This is equivalent to
\[
|A|(1-2\ep)>\Sigma(1+2\ep+\de).
\]
Plugging in Inequality (\ref{ineq:tiling}), we see that it is enough to prove that
\[
(1-2\ep)>(1+2\ep+\de)\frac{1-2\de}{1-\de},
\]
which is our assumption on $\ep$.
\end{proof}

By the lemma, we see that if $|A_k|>2\de|A|$, then $|F_k\tB_k\cap A|\geq\frac{\de}{4c}|B_k|\geq\frac{\de}{8c}|A_k|$, which implies that $|A_{k-1}|<(1-\frac{\de}{8c})|A_{k})$, and, using induction, that $|A_{k-1}|<(1-\frac{\de}{8c})^{N-k}|A|$. By our assumption on $N$, we get that this cannot hold for all $k$. Hence $|A_1|<2\de|A|$, which shows the third conclusion of the proposition.
\end{proof}

\subsection{Ergodic Theorem for Castles}

We fix a cross-section $(S,\cB ,\mu ,\cR ,\al )$.

\begin{defn} Let $A\subset S$ be a Borel subset. A castle with base $A$ is a Borel subset $\cT\subset\cR$ such that
\begin{enumerate}
\item The set $\{(x,x)|x\in A\}$ is contained in $\cT$.
\item For all $x,y\in S$, if $(x,y)\in\cT$ then $x\in A$.
\item For every $x\in A$, the set $\cT_x:=\{y\in S\q|\q(x,y)\in\cT\}$ is finite.
\item For every two distinct points $x_1,x_2\in A$, the sets $\cT_{x_1}$ and $\cT_{x_2}$ are disjoint.
\end{enumerate}
We call $\cT_x$ the {\em tower over $x$}. Note that if $\cT$ is a castle with base $A$, then the function $x\mapsto\cT_x$ from $A$ to the finite subsets of $S$ is Borel measurable. The {\em range} of the castle is the the set $rg(\cT)=\cup_{x\in A}\cT_x$. If $\cG$ is a sub-sigma-algebra, we say that the castle $\cT$ is $\cG$-measurable if $\cT$ is $\cG$-measurable (which implies that the base of $\cT$ is $\cG$-measurable). Finally, we say that $\cT$ covers $\de$ of $S$, if $\mu(rg(\cT))>\de$.
\end{defn}

\begin{rem} In the Ornstein--Weiss theory, there are several notions (casle, real-tower, $\ep$-quasi tower) similar to our notion of castle. The definition given here is more restrictive than $\ep$-quasi tower in that we require that the sets $\cT_x$ are disjoint, but less restrictive in that we do not impose conditions on the possible values of $\al(\cT)$.
\end{rem}

\begin{defn} Let $\cT$ be a castle. Define a measure $\mu_\cT$ on its base by
\[
\mu _{\cT} (B)=\int _B |\cT _x |d\mu (x).
\]
In words, $\mu_\cT(B)$ is the $\mu$-measure of the part of the range of the castle that lies over $B$.
\end{defn}

\begin{defn} Given $T\subset S$ and a compact $K\subset G$, the $K$-interior of $T$ is the set 
\[
int_K T=\{ t\in T | (\forall r \in S) ( \al (t,r)\in K \longrightarrow r\in T) \} .
\]
The $K$-boundary of $T$ is $\partial _K T=T\setminus int_K T$.
\end{defn}

\begin{defn} Let $\cT$ be a castle with base $A$ and let $K\subset G$ be compact set. For $\ep >0$ we say that $\cT$ is $(K,\ep )$-invariant if 
\[
\mu _{\cT} \left ( \left \{ s\in A \left | \frac{|\partial _K \cT _s |}{|\cT _s |} >\ep \right \}  \right ) <\ep \mu_{\cT}(A)\right . .
\]
\end{defn}

\begin{thm} \label{erg:thm:castles} (Ergodic theorem for castles.) Let $(S,\cB ,\mu ,\cR ,\al )$ be a cross-section. Let $h:S\to \bC$ be a bounded function. For every $\de >0$ there is an $\eta >0$ and a compact set $K\subset G$ such that if $\cT$ is a $(K,\eta)$-invariant castle with base $A$ such that $\mu (rg(\cT ))> \de$, then there is a subset $B\subset A$ such that $\mu _{\cT} (B)>(1-\de )\mu _{\cT} (A)$ and such that for any $s\in B$,
\[
\left | \frac{1}{|\cT _s|}\sum _{t\in \cT _s} h(t) -\int _S h d\mu \right | <\de .
\]
\end{thm}

\begin{proof} Assume that $h$ is bounded by $M$ and that it has zero mean. Given $\de >0$, let $N$ and $\ep$ be as in Proposition \ref{prop:tiling}. Given a sequence of (F\o lner) sets $F_1 ,\dots ,F_N$, a castle $\cT$, and a point $s$ in the base of $\cT$, let $B_s$ be the set of points $t\in \cT _s$ such that for all $i$, both $t\in int_{F_i} \cT _s$ and 
\begin{equation} \label{eq:erg:thm:castles:mean}
\left | \frac{1}{|\{ u | \al (t,u)\in F_i\} |}\sum _{ u : \al (t,u)\in F_i} h(u) \right | <\de .
\end{equation}
We call $s$ good if $(\al (s,\cT _s),\al (s,B_s), F_i)$ satisfies the conditions of the tiling lemma. Suppose $s$ is good. Then we can find $t_1,\dots ,t_K \in \cT _s$ and a function $n:\{1,\dots ,K\}\to\{1,\dots ,N\}$ such that the sets $C_i=\{ u\in \cT _s | \al(s,u) \in F_{n(i)}\al(s,t_i)\}$ are $\de$-disjoint and cover $1-\de$ of $\cT _s$. Let $D_i=C_i \setminus \cup _{j<i}C_j$ and let $D=\cup D_i=\cup C_i$. Then the $D_i$ are disjoint, $|D_i|>(1-\de)|C_i|$, and $|D|/|\cT _s|>1-\de$. In this case,
\[
\left | \frac{1}{|\cT _s |}\sum _{t\in \cT _s} h(t) \right | \leq \frac{1}{|\cT _s|}\sum _{t\in \cT _s \setminus D} |h(t)|+\frac{1}{|\cT _s|}\sum _i \left | \sum _{t\in D_i} h(t) \right | \leq
\]
\[
M\de +\frac{1}{|\cT _s|}\sum _i \left ( \left | \sum _{t\in C_i}h(t) \right | +\sum _{t\in C_i\setminus D_i} |h(t)| \right ) < (2M+1)\de .
\]

We now show that we can choose F\o lner sequence $F_i$, a compact set $K\subset G$, and $\eta >0$ such that if $\cT$ is a $(K,\eta)$-invariant castle, then most (in the sense of $\mu _{\cT}$-measure) of the points of the base of $\cT$ are good. First note that for every compact set $U$, if $F$ is $(U,1)$-invariant and $X\subset F$ is the set of centers of disjoint right translates of $U$, then
\[
|X|\la(U)=\la(XU)\leq\la(FU)\leq2\la(F),
\]
and so $|X|\leq\frac{2}{\la(U)}\la(F)$.

It is therefore enough to show that for most points $s$ in the base of $\cT$ and for most points $t\in \cT _s$,
\begin{enumerate}
\item $\{ u | \al (t,u)\in F_i \} \subset \cT _s$ for all $i$.
\item \eqref{eq:erg:thm:castles:mean} holds.
\item $|\{ u\in \cT _s | \al (t,u)\in (\cup _{j<i} F_j)^{-1} F_i \} |<(1+\ep)|\{ u\in \cT _s | \al (t,u)\in F_i \}|.$
\end{enumerate}

Choose the $F_i$ to be sufficiently invariant so that the mean ergodic theorem is satisfied for $1-\de \ep ^2$ of the points in $S$, for each $F_i$, and for the functions $h$ and $1$. By Fubini, 2. is satisfied for at least $1-\ep ^2$ of the points in $rg(\cT )$. By Fubini again, for $1-\ep$ of the points of the base, 2. is satisfied for $1-\ep$ of the points above them. We require additionally that $F_i$ is $((\cup _{j<i} F_j)^{-1} ,\frac{\la(U)}{2}\ep )$-invariant, and denote $L_i = (\cup _{j<i} F_j)^{-1} F_i \setminus F_i$. We can further require that $L_i$ is $(U,1)$-invariant, and so for all $t\in \cT _s$, $|\{ u \in \cT _s | \al (t,u)\in L_i \} |<c\la (L_i)$, where $c=\frac{2}{\la(U)}$. By the ergodic theorem for the function $1$, $|\{ u | \al (t,u) \in F_i \} |>(1-\ep )\la (F_i )$. Now we add the demand that $\cT$ is $(\cup _i F_i , \ep )$-invariant. Then 1. is satisfied for most pairs $(s,t)$, and for each point $t\in int_{F_i}\cT _s$,
\[
\frac{|\{ u \in \cT _s | \al (t,u) \in L_i \} |}{|\{ u\in \cT _s | \al (t,u)\in F_i \} |} < \frac{c\la (L_i)}{(1-\ep )\la (F_i)} < \ep
\]
which implies 3.
\end{proof}

\section{Entropy Theory} \label{sec:ent.theory}

\subsection{Entropy Theory for Amenable Groups Actions}

The entropy theory of actions of amenable groups is developed in \cite{OW}. The definition of entropy takes simple form for groups with zero self entropy, which are defined as follows:
\begin{defn} Let $G$ be a unimodular amenable group, and choose a right invariant metric, $d$, on $G$ that generates the topology. We say that $G$ has zero self entropy if the following condition holds: For every compact subset $K$ of $G$, and for every $\ep >0$, if $F\subset G$ is sufficiently invariant, then there is a partition of $K$ into less than $2^{\ep \la (F)}$ sets such that for every two points $g_1 ,g_2$ in the same part and every $f\in F$, we have $d(fg_1 ,fg_2)<\ep$.
\end{defn}

\begin{exa} \begin{enumerate}
\item Every nilpotent Lie group has zero self entropy.
\item There are unimodular solvable Lie groups that do not have zero self entropy.
\end{enumerate}
\end{exa}

We now briefly remind the reader the definition of entropy for actions of groups with zero self entropy. We are given an action of the group $G$ on a probability space $X$. A partition of $X$ is a function $P$ from $X$ to a finite set, whose size we denote by $|P|$. For every partition $P$ of $X$ and a compact set $F\subset G$, we say that $x,y\in X$ are $(P,F,\ep )$-close if 
\[
\la ( \{ f\in F | P(fx)\neq P(fy)\} )< \ep \la (F).
\]
A $(P,F,\ep )$-ball is a subset $A\subset X$ such that every two points in $A$ are $(P,F,\ep )$-close. The entropy of the action with respect to $P$ is defined to be the minimal number $h=h(G,P)$ such that for every $\ep>0$, if $F\subset G$ is sufficiently invariant, then there is a collection $\cC$ of $(P,F,\ep)$-balls that covers $1-\ep$ of $X$ and has size less than $2^{(h+\ep )\la (F)}$. We define the entropy of the action of $G$ on $X$, denoted by $h(X)$, to be the supremum of the numbers $h(X,P)$, where $P$ is taken from the set of partitions of $X$.

\begin{rem} \begin{enumerate}
\item If $X$ is a Borel space endowed with regular probability measure, then $h(X)$ is the supremum of the numbers $h(X,Q)$, where $Q$ is taken from the set of partitions of $X$ such that any part of $Q$ has negligible boundary. This is true since this collection of partitions generates the sigma-algebra.
\item In the other direction, one can consider more general partitions, which are (Borel) maps to some compact metric space $(Z,d)$. If $P$ is such a generalized partition, we say that two points $x,y\in X$ are $(P,F,\ep)$-close if
\[
\int _F d(P(fx),P(fy)) < \ep \la (F),
\]
and define similarly the notions of $(P,F,\ep)$-ball and the entropy of the action with respect to $P$. It turns out (see \cite{OW}) that the supremum of $h(X,P)$ where $P$ is taken from the set of generalized partitions is again $h(X)$.
\end{enumerate}
\end{rem}

\subsection{Entropy Theory for Cross Sections}

In our definition of entropy for cross-sections, we choose a slightly different route, concentrating on castles.

\begin{defn} \label{defn:T.ball} Let $(S,\cB,\mu,\cR,\al)$ be a cross-section and let $\cG$ be a sub-sigma-algebra of $\cB$. Let $P$ be a finite partition of $S$, not necessarily $\cG$-measurable, let $\cT$ be a castle with base $A$, and let $\ep >0$. A {\em $(P,\cT ,\ep ,\cG)$-ball} is a triple $(B,E,\phi )$ such that $B$ is a subset of $A$, not necessarily $\cG$-measurable, $E$ is a finite set, and $\phi :B\times E \to S$ is a restriction of a $\cG$-measurable, one-to-one function such that
\begin{enumerate}

\item $\phi (x,e)\in \cT _x$ for all $x\in B$, $e\in E$.

\item $|\phi (x,E)|>(1-\ep )|\cT _x |$ for every $x\in B$.

\item $x\mapsto P(\phi (x,e))$ is constant for every $e\in E$.

\end{enumerate}
\end{defn}

Given $(B,E,\phi )$ as above, for any $s,t\in B$ we get a partially defined map $\phi _{s,t}$ between $\cT _s$ and $\cT _t$ by $\phi _{s,t}(r)=\phi (t,\pi _2 (\phi ^{-1}(r)))$ where $\pi _2$ is the projection to the second coordinate. The function $\phi _{s,t}$ is defined for at least $1-\ep$ of the points of $\cT _s$, and its image consists of at least $1-\ep$ of the points of $\cT _t$. Conversely, if $B$ is given and we can show that for every $s,t \in B$ there is such a $\phi _{s,t}$ that depends on $s$ and $t$ in a $\cG$-measurable way, then $B$ can be decomposed in a $\cG$-measurable way into balls---that is, there are $\cG$-measurable sets $C_1,\dots ,C_N$, such that $B\subset \cup C_i$, and there are $E_i ,\phi _i$ as above such that $(B\cap C_i, E_i ,\phi _i )$ are $(P,\cT ,2\ep ,P)$-balls.

\begin{defn} Let $\cG$ be a sigma-algebra, let $A\subset S$ be a $\cG$-measurable set, and let $\cC =\{ C_1,\dots ,C_N\}$ be a collection of Borel sets (not necessarily $\cG$-measurable) contained in $A$. Denote $\log _{+} (x)=\max \{\log(x),0\}$. The relative logarithmic size of $\cC$ with respect to $\cG$ is
\[
rls_A(\cC |\cG )=\int _A \log _+(| \{ i|E(1_{C_i}|\cG )(x)>0\} |)d\mu (x).
\]
\end{defn}

\begin{rem} \label{rem:G.fibers} One can change the sets $C_i$ by sets of measure 0 such that the integrand $| \{ i|E(1_{C_i}|\cG )(x)>0\} |$ is the number of sets in $\cC$ that intersect the $\cG$-fiber of $x$. 
%Conversely, the following will be usefull: suppose that for each $x\in A$, there is a collection $\cC _x$ of subsets of the fiber of $x$ with less than $H$ elements. Assume also that the map $x\mapsto \cC _x$ is $\cG$ measurable. Then there is a finite collection $\cD =\{ D_1 ,\dots ,D_N\}$ of sets with relative logarithmic size $\mu (A)\log H$ such that for almost any $x\in A$, the non empty sets in the intersection of $D_i$ and the fiber of $x$ are the sets in $\cC _x$.
\end{rem}

\begin{defn} Suppose that $\cT$ is a $\cG$-measurable castle with base $A$. Let $P$ be a partition. Denote by $h_{\cT}^{\ep} (P|\cG )$ the infimum of $rls_A (\cC )$ where $\cC$ is a collection of $(P,\cT ,\ep ,\cG)$-balls that cover $A$.
\end{defn}

\begin{rem} Clearly, $h_{\cT} ^{\ep}$ is monotone in $\cT$---if $\cT$ and $\cT '$ are $\cG$-measurable castles with disjoint ranges, then for every partition $P$ and $\ep >0$, we have $h_{\cT \sqcup \cT '} ^{\ep} (P|\cG )=h_{\cT} ^{\ep} (P|\cG )+h_{\cT '} ^{\ep} (P|\cG )$. Also, if $|P|$ is a set whose size is less than $N$, then $h_{\cT}^{\ep}(P|\cG)<\mu (rg(\cT ))\log N$.
\end{rem}

\begin{thm} \label{comp:ent:castles} Let $(S,\cB ,\mu ,\cR ,\al )$ be a cross-section. Suppose that $\cG$ is a sigma-algebra such that $\al$ is $\cG$-measurable. Let $0<\ep<1$. Suppose that $\cT$ is a tower that covers more than $1-\ep$ of $S$. There is a compact set $L\subset G$ and a $\de >0$ such that if $\cT '$ is an $(L,\de )$-invariant castle that is $\cG$-measurable and covers $1-\ep$ of $S$, then $h_{\cT '}^{4\ep }(P|\cG )\leq h_{\cT}^{\ep}(P|\cG )+2\ep+\ep \log (|P|)$.
\end{thm}

\begin{proof} Let $A$ be the base of $\cT$. There is a compact subset $K\subset G$ such that for all $s$ in $A$, except for a subset of $\mu _{\cT}$-measure as small as we want, $\al (s,\cT _s )\subset K$. After removing this small set (and noting that we can choose $K$ to be so large so that we still have a castle that covers more than $1-\ep$ of $S$), we may assume that $\al (s,\cT _s )\subset K$ for all $s$ in $A$. Let $\cC$ be a collection of $(P,\cT,\ep,\cG)$-balls that covers $A$, such that $rls_A(\cC|\cG )<h_{\cT}(P|\cG )+\ep$. Let $H(x)=\log_+ (|\{ i | E(1_{C_i}|\cG )(x)>0\} |)$ (so $H(x)=0$ outside $A$). By the ergodic theorem for castles, if $\cT'$ is a sufficiently invariant castle with base $A'$ that covers more than $\ep$ of $S$, then for a subset $B\subset A'$ of $\mu _{\cT '}$-measure greater than $1-\ep$, we have that for every $s\in B$,
\begin{equation} \label{eq:ent:castles:H}
\left | \frac{1}{|\cT _s '|}\sum _{t\in \cT_s '}H(t)-\int_S Hd\mu \right | < \ep
\end{equation}
and
\begin{equation} \label{eq:ent:castles:1}
\left | \frac{1}{|\cT _s'|}\sum _{t\in \cT _s '}1_{rg(\cT )}(t) - \mu (rg (\cT )) \right | <\ep .
\end{equation}
Taking $\cT'$ as sufficiently invariant, and by making $B$ smaller, we can assume that for all points $s$ of $B$, both (\ref{eq:ent:castles:H}) and (\ref{eq:ent:castles:1}) hold, that $\cT'_s$ is $(K,\ep)$-invariant, and that $\mu_{\cT'}(B)>(1-\ep)\mu_{\cT'}(A)$. 

Since $\cT '$ and $B$ are $\cG$-measurable, it is possible to divide $B$ into $\cG$-measurable subsets, $B_i$, and to find a one-to-one and measure-preserving map $\phi _i :B_i \times E_i \to \cT '$, where $E_i$ are finite sets, such that $\phi _i (s,E_i) = \cT '_s$ and the following holds: if we denote for every $s,t \in B_i$ the identification $t\mapsto \phi _i (s_2, \pi _2 (\phi _i ^{-1}(t)))$ between $\cT '_s$ and $\cT '_t$ by $\phi _{s,t}$, then $E(1_{C_i}|\cG )(r) >0$ iff $E(1_{C_i}|(\cG )(\phi _{s,t}(r)) >0$ for all $s,t\in B_i$ and $r\in \cT '_s$. Let $s\sim t$ if, for every $r\in \cT _s '$, the points $r$ and $\phi _{s,t}(r)$  are contained in the same $C_i$. Denote the equivalence classes inside $B_j$ by $\cD ^j _1 ,\dots ,\cD ^j _{N_j}$. If $s\in B_j$, the number $N_j$ of equivalence classes is $\exp (\sum _{t\in \cT _s} H(t))$, and by \eqref{eq:ent:castles:H}, this number is less than $\exp ((h_\cT^\ep(P|\cG)+2\ep)|\cT'_s|)$. Therefore, the relative logarithmic size of of the collection $\{\cD^j_i\}$ over $B$ is less than $h_\cT(P|\cG)+2\ep$.

Next, we show that each $\cD _i ^j$ is a $(P,\cT ',4\ep,\cG)$-ball. For $s\in D_i ^j$, denote $R_s=\{ r\in \cT '_s | r\in A , \cT _r \subset \cT '_s\}$. From \eqref{eq:ent:castles:1} and the assumption that $\cT'_s$ is $(K,\ep)$-invariant, we get that for all $s\in B$,
\[
\frac{|\bigcup _{r\in R_s} \cT _r |}{|\cT '_s|} > 1-3\ep.
\]

Given $s_1,s_2$ in $B$, we define a map $\psi : \cT '_{s_1} \to \cT '_{s_2}$. On $R_{s_1}$ it is equal to $\phi _{s_1,s_2}$. Given $r\in R_{s_1}$, since $r$ and $\phi _{s_1,s_2}(r)$ are in the same $C_i$, there exists a function $\theta _r:\cT _r \to \cT _{\phi _{s_1,s_2}(r)}$ such that $P(t)=P(\theta _r (t))$ for $1-\ep$ of the points in $\cT _r$. Since the sets $\cT _r$ for $r\in R_s$ are disjoint, we can define a function $\theta : \cup _{r\in R_{s_1}}\cT _r \to \cT '_{s_2}$ such that $P(\theta (t))=P(t)$ for $1-\ep$ of the points $t\in \cup _{r\in R_{s_1}}\cT _r$. Extending this function in an arbitrary way to the rest of $\cT '_{s_1}$, we get a $\cG$-measurable function such that $P(w)=P(\theta (w))$ for $1-4\ep$ of the points.

Therefore, if we denote the restriction of $\cT'$ to $B$ by $\cT'|_B$, we get that
\[
h^{4\ep}_{\cT'}(P|\cG) \leq h^{4\ep}_{\cT'|_B}(P|\cG)+\mu(rg(\cT'_{A'\setminus B}))\log|P|\leq h_\cT(P|\cG)+2\ep+\ep\log|P|.
\]
\end{proof}

\begin{defn} Let $S$ be a cross-section, let $P$ be a partition and let $\cG$ be a sub-sigma-algebra. We define the relative entropy to be
\[
h(P|\cG ) = \lim _{\ep \to 0} \lim _{\cT} h_{\cT} ^{\ep}(P|\cG )
\]
where the inner limit is taken over more and more invariant castles that cover $1-\ep$ of $S$.
\end{defn}

\begin{rem} By the previous theorem, the limit above exists.
\end{rem}

The following is an analogue of Shannon--McMillan's theorem.

\begin{thm} \label{thm:SM} Suppose $(S,\cB ,\mu ,\cR ,\al )$ is a $U$-discrete cross-section. Let $P$ be a partition of $S$, and let $\cG\subset\cB$ be a sub-sigma-algebra such that $\al$ is $\cG$-measurable. Then $h(P|\cG )$ is the infimum over the positive real numbers $h$, such that for every $\ep>0$, for sufficiently invariant F\o lner sets $F\subset G$, there is a collection $\cC =\{ C_i ,\dots ,C_N\}$ of subsets of $S$, finite sets $E_1 , \dots ,E_N$, and one-to-one functions $\phi_i:C_i \times E_i \to S$, that are restrictions of $\cG$-measurable functions, such that
\begin{enumerate}
\item For each $i$ and $s\in C_i$, we have that $\phi_i(s,E_i)\subset \{ t | \al (s,t)\in F\}$ and $|\phi (s,E_i )|>(1-\ep )|\{ t | \al (s,t)\in F\} |$.
\item For each $i$ and $e\in E_i$, the function $s\mapsto P(\phi_i(s,e))$ ($s\in C_i$) is constant.
\item $\left|\{i|E(1_{C_i}|\cG)(s)>0\}\right|<\exp((h+\ep)|\{t|\al(s,t)\in F\}|)$ for all $s\in S$.
\item $\mu \left ( \bigcup _i C_i \right ) >1-\ep .$
\end{enumerate}

\end{thm}

\begin{proof} Choose a neighborhood $V$ such that $V^2\subset U$. In one direction, assume that $h$ is such that for any $\ep>0$, if $F$ is sufficiently invariant, then there are $C_i,E_i,\phi_i$ as above. We show that for any $\ep>0$, if $\cT$ is a sufficiently invariant castle that covers $1-\ep$ of $S$, then $h_\cT^\ep(P|\cG)<h+\ep$. 

By Proposition \ref{prop:tiling}, there is $\eta>0$ and $n>0$ such that if $A,B,F_1,\ldots,F_n\subset G$ satisfy the conditions of the tiling lemma with parameter $\eta$, then the conclusions of Proposition \ref{prop:tiling} hold with $\de$ replaced by $\ep$. Choose $F_1,\ldots,F_n$ such that
\begin{itemize}
\item $F_i$ is $(U,1)$-invariant.
\item $\la(F_i)\geq0$.
\item For every $i$, $(\cup_{j<i}F_j)^{-1}F_i\setminus F_i$ is $(U,1)$-invariant.
\item For every $i$, $\la((\cup_{j<i}F_j)^{-1}F_i\setminus F_i)<\eta\la(F_i)$.
\end{itemize}
 For every $1\leq k\leq n$, let $\cC^k=\{C^k_1,\ldots,C^k_{m(k)}\}$ be a collection of subsets of $S$, let $E^k_i$, $1\leq i\leq m(k)$, be a collection of finite sets, and let $\phi^k_i:C^k_i\times E^k_i\to S$, $1\leq i\leq m(k)$, be a collection of one-to-one functions that are the restrictions of $\cG$-measurable functions, such that {\em 1.,2.,3.,} and {\em 4.} in the conditions of the theorem are satisfied with $\ep$ replaced by $\ep\de/n$.

By the ergodic theorem, if $\cT$ is sufficiently invariant, then for $1-\ep$ of the points $s$ in the base of $\cT$ we have that at least a $1-\eta$ portion of $\cT_s$ is contained in $\cap_{k=1}^n\cup \cC^k$. For each such a point, apply Proposition \ref{prop:tiling} to find
\begin{itemize}
\item A number $K(s)$.
\item Points $x_i(s)\in\cT_s\cap\left(\cap_k\cup\cC^k\right)$, for $1\leq i\leq K(s)$.
\item Indices $j(i,s)\in\{1,\ldots,n\}$.
\end{itemize}
such that the sets $\{y|\al(x_i(s),y)\in F_{j(i,s)}\}$ are $\ep$-disjoint and cover $1-\ep$ of $\cT_s$. From the proof of Proposition \ref{prop:tiling}, we can assume that the functions $K(-),x_i(-),$ and $j(i,-)$ are $\cG$-measurable. For two points, $s$ and $t$, in the base of $\cT$ for which the above holds, we say that $s$ is equivalent to $t$ if $K(s)=K(t)$, for all $1\leq i\leq K(s)$, $j(i,s)=j(i,t)$, and for all $1\leq i\leq K(s)$, the points $x_i(s)$ and $x_i(t)$ belong to the same set in $\cC^{j(i,s)}$. It is easy to see that each equivalence class is a $(P,\cT,2\ep,\cG)$-ball. Moreover, for any $s$, the number of equivalence classes that intersect the $\cG$-fiber of $s$ is at most
\[
\prod_{i=1}^{K(s)}\left|\{k|E(1_{C_k^{j(i,s)}}|\cG)(x_i(s))>0\}\right| \leq
\]
\[
\leq \exp\left((h+\ep)\sum_{i=1}^{K(s)}\left|\{t|\al(x_i(s),t)\in F_{j(i,s)}\}\right|\right) \leq
\]
\[
\leq \exp\left( (h+\ep)(|\cT_s|+\ep)\right) \leq \exp((h+2\ep)|\cT_s|),
\]
and therefore $h_\cT^{2\ep}(P|\cG)\leq h+2\ep$.

In the other direction, let $\cT$ be a $\cG$-measurable castle with base $A$ that covers $1-\ep$ of $S$. After discarding a small enough subset of the base of $\cT$, we can assume that the set $\{\al(s,t)|(s,t)\in\cT\}$ is contained in some compact set $L$. Suppose that there is a collection $\cD=\{D_1,\ldots,D_n\}$ of $(P,\cT,\ep,\cG)$-balls that covers $1-\ep$ of $A$, such that $rls_A\cD<h$. Let $H:S\to\bR$ be the function such that $H(s)$ is equal to $\log_+(|\{i|E(1_{D_i}|\cG)(s)>0\}|)$ if $s$ belongs to $A$, and is equal to $0$ if $s$ is not in $A$. We have that $\int_SHd\mu=rls_A\cD<h$. Let $F\subset G$ be sufficiently invariant, such that for $1-\ep$ of the points $s$ of $S$, 
\begin{itemize}
\item $(1-\ep)\la(F)\leq|\{t|\al(s,t)\in F\}|\leq(1+\ep)\la(F)$.
\item $|\{t\in S | \al(s,t)\in F \wedge t\in rg(\cT)\}|>(1-\ep)|\{t\in S|\al(s,t)\in F\}|$.
\item $\sum_{t : \al(s,t)\in F}H(t)\leq (h+\ep)|\{t|\al(s,t)\in F\}|$.
\end{itemize}
The existence of such $F$ is given by the ergodic theorem.

By enlarging $F$, we can assume that $F$ is $(L,\ep)$-invariant. For all but $1-2\ep$ of the points $s$ of $S$, there is a number $K(s)$ and elements $x_i(s)\in S$, $1\leq i\leq K(s)$ such that
\begin{itemize}
\item For every $i$, the point $x_i(s)$ belongs to the base of $\cT$.
\item The sets $\cT_{x_i(s)}$ cover $1-2\ep$ of $\{t|\al(s,t)\in F\}$.
\item $\sum_{i=1}^{K(s)}H(x_i(s))\leq (h+2\ep)|\{t|\al(s,t)\in F\}|$.
\end{itemize}
We can assume, moreover, that the functions $K(-)$ and $x_i(-)$ are $\cG$-measurable. Define an equivalence relation on the set of points of $S$ for which the above holds: $s\sim t$ if $K(s)=K(t)$ and for all $1\leq i\leq K(s)$, the points $x_i(s)$ and $x_i(t)$ belong to the same set in $\cD$. 

Let $\cC=\{C_i\}$ be the collection of equivalence classes. For every $K$, the number of equivalence classes that are contained in $\{s|K(s)=K\}$ is less than $2^{(h+3\ep)\la(F)}$. Therefore, for every $s\in S$,
\[
|\{i|E(1_{C_i}|\cG)(s)>0\}|\leq2^{(h+3\ep)\la(F)}\leq2^{(h+4\ep)|\{t|\al(s,t)\in F\}|}.
\]

For every two equivalent points $s,t$, there is a map $\phi_{s,t}:\{u|\al(s,u)\in F\}\to\{t|\al(t,v)\in F\}$, such that, for $1-3\ep$ of the points in $\{u|\al(s,u)\in F\}$, we have $P(\phi_{s,t}(u))=P(u)$. As explained after Definition \ref{defn:T.ball}, there is a $\cG$-measurable refinement of the $C_i$'s, such that, for each new part $Y$, there is a set $E$ and a function $\phi:X\times E\to S$, which is the restriction of a $\cG$-measurable function, such that 1. and 2. in the theorem hold. Since the refinement is $\cG$-measurable, the relative size of it is unchanged, and so 3. holds. Lastly, 4. holds because the union is just the union of $\cC$.

\end{proof}

\begin{thm} \label{transference} Let $\cR$ be a probability-preserving equivalence relation on the probability space $(S,\cB ,\mu )$. Let $\cG$ be a sub-sigma-algebra of $\cB$, and suppose that there are two amenable groups $G,H$, with zero self entropy, and two $\cG$-measurable cocycles $\al :\cR \to G$ and $\be :\cR \to H$ such that $(S,\cB ,\mu ,\cR ,\al )$ and $(S,\cB ,\mu ,\cR ,\be )$ satisfy the mean ergodic theorem. For every partition $P$, denote by $h_{\al}(P|\cG )$ the relative entropy of $P$ with respect to $\cG$ for the cross-section $(S,\cB ,\mu ,\cR ,\al )$, and define similarly $h_{\be}(P|\cG )$. Then for each partition $P$, $h_{\al}(P|\cG )=h_{\be}(P|\cG )$.
\end{thm}

\begin{proof} Let $S,\cB,\mu,\cR,\al,\be,\cG$ be as in the statement of the theorem.
\begin{lem} For every $\ep>0$ and a compact set $L\subset H$, there is a compact set $F\subset G$ such that if $\cT$ is a castle, and is $(F,\ep)$-invariant, then it is also $(L,2\ep)$-invariant.
\end{lem}
\begin{proof} Let $K\subset G$ be a compact set such that the set
\[
Z=\{x|\textrm{for every $y$ such that $\be(x,y)\in L$, we have $\al(x,y)\in F$}\}
\]
has measure greater than $1-\ep/2$. By the ergodic theorem, there is a compact set $F\subset G$ that contains $K$ such that if $\cT$ is $(F,\ep)$-invariant, then $1-\ep$ of the points in the range of $\cT$ are in $Z$. But if $s$ is in the base of $\cT$ and $x\in\cT_s$ is both in $Z$ and in the $F$-interior of $\cT_s$, then $x$ is also in the $L$-interior of $\cT_s$.
\end{proof}

Let $L_n\subset H$ be a F\o lner sequence in $H$. Take a F\o lner sequence $F_n\subset G$ that satisfies the conclusion of the lemma.

The definition of a $(P,\cT,\ep,\cG)$-ball does not depend on the cocycle, and hence neither is the definition of $h_\cT^\ep(P|\cG)$. The number $h_\al(P|\cG)$ is the limit of $h_\cT^\ep(P|\cG)$ as $\cT$ is taken to be $(L_n,\ep)$-invariant, and $h_\be(P|\cG)$ is the limit of $h_\cT^{2\ep}(P|\cG)$ as $\cT$ is taken to be $(F_n,2\ep)$-invariant. These limits are equal.
\end{proof}

The following is a version of Abramov's theorem.
\begin{thm} \label{Abramov} Let $G$ be an amenable group of zero self entropy. Assume $G$ acts on a Borel probability space $(X,\cB,m)$, and that $S\subset X$ is a $U$-discrete cross-section, for some neighborhood $U$ of $1$ in $G$. Let $\cG$ is a $G$-invariant sub-sigma-algebra such that $S$ is $\cG$-measurable. Then
\[
h(G,X|\cG )=h(S|\cG _{S})
\]
where $\cG _S$ denotes the restriction of $\cG$ to $S$.
\end{thm}

\begin{rem} \begin{enumerate}
\item On the left hand side, the relative entropy $h(G,X|\cG)$ is defined as the difference of the entropy of the action of $G$ on $(X,\cB)$ and the entropy of the action of $G$ on $(X,\cG)$. Since $G$ has zero self entropy, the relative entropy is the difference of the spatial entropies (see \cite[II.4]{OW}). For our argument, we need that the relative entropy can be approximated by choosing a large Folner set $L$, and a $\cG$-measurable partition of $X$, such that most parts can be partitioned into $2^{(1+\ep)\la(L)h(G,P|\cG)}$ $(P,L,\ep)$-balls.
\item The formula is simpler than in the classical Abramov theorem because of our normalization \ref{intensity1} of the Haar measure.
\end{enumerate}
\end{rem}

\begin{proof} Let $P$ be a partition of $X$, and let $\ep>0$. We assume that every part of $P$ has negligible boundary. Choose a compact set $F\subset G$ such that $m(F\cdot S)>1-\ep/2$. Let $\eta=\frac{\ep}{2\la(F)}$. Define a generalized partition $\tP$ as follows: the value set of the partition is the set of functions from $F$ to $|P|$, with metric given by
\[
d(\phi,\psi)=\la\{f\in F | \phi(f)\neq\psi(f)\},
\]
and $\tP (x)$ is the function $f \mapsto P(fx)$. Let $\Om$ be the image of $\tP$. Since $\tP$ is continuous and $X$ is compact, $\Om$ is also compact. Choose an $\eta/2$-net $\cA$ in $\Om$. Choose also a measurable map $\rho:\Om\to\cA$ such that $d(\om,\rho(\om))<\eta/2$ for all $\om\in\Om$. Let $\ttP$ be the partition of $S$ given by $\ttP(s)=\rho(\tP(s))$.

Let $W$ be a neighborhood of the identity in $G$ such that the set
\[
B=\{x\in X | \textrm{the function $w\mapsto P(wx)$ is constant on $w\in W$}\}
\]
has measure greater than $1-\ep/2$.

Finally, choose $L\subset G$ compact such that $L$ is $(F,\ep)$-invariant, and for all $s\in S$, outside a set of $\mu$-measure less than $\ep$, the following hold:
\begin{enumerate}
\item $\la\{f\in L| fs\in B\}>(1-\ep)\la(L)$.
\item \label{asump:FS} $\la\{f\in L| fs\not\in F\cdot S\}<\ep\la(L)$.
\item $|\{f\in L|fs\in S\}|<2\la(L)$.
\item (Theorem \ref{thm:SM}) There is a collection $\cC=\{C_i\}$ of subsets of $S$ such that $rls(\cC|\cG)<\la(L)(h(S|\cG)+\ep)$, such that $\mu(\cup\cC)>1-\ep$, and such that each $C_i$ is a $(\ttP,L,\eta)$-ball.
\item (Zero self-entropy of $G$) There is a partition of $F$ into less than $2^{\ep \la (L)}$ parts $M_i$ such that for $g_1,g_2$ in the same $M_i$ and for any $f\in L$, there is $\de (f,g_1 ,g_2)\in W$ such that $fg_1=\de (f,g_1 ,g_2)fg_2$.
\end{enumerate}

\begin{lem} There is a $\cG$-measurable partition $\cD$ such that every set in the common refinement $\cC\vee\cD$ is a $(P,L,3\ep)$-ball.
\end{lem}

\begin{proof} Suppose first that $x_1,x_2\in C_i$ have the same return times to $S$, i.e.
\[
A:=\{f\in L|fx_1\in S\}=\{f\in L|fx_2\in S\}.
\]
By assumption, there is a set $B\subset A$ such that $|B|>(1-\eta)|A|$ and for every $f\in B$, $\ttP(fx_1)=\ttP(fx_2)$. If we denote by $1_\phi$ the characteristic function of $\phi$, then
\[
\la\{f\in L|P(fx_1)\neq P(fx_2)\}=\int_L 1_{P(fx_1)\neq P(fx_2)}df \leq
\]
\[
\la(L\setminus FA)+\la(FA\setminus FB)+\sum_{b\in B}\int_{F}1_{P(fbx_1)\neq P(fbx_2)}df.
\]
The first summand is less than $\ep\la(L)$ by assumption \ref{asump:FS} on $L$. The second summand is less than
\[
|A\setminus B|\la(F)<\eta|A|\la(F)\leq\ep\la(L),
\]
and since for every $b\in B$, the distance between the $(\tP,F)$-names of $bx_1$ and $bx_2$ is less than $\eta$, we get that the third summand is less than $|B|\eta\leq\ep\la(L)/2$. Therefore,
\[
\la\{f\in L|P(fx_1)\neq P(fx_2)\}<2.5\ep\la(L).
\]

Since the cocycle $\al$ is $\cG$-measurable, for every small enough neighborhood $Z\subset G$ of the identity, there is a $\cG$-measurable partition $\cD_Z$ such that every $x_1,x_2$ in the same part and every $f\in L$ such that $fx_1\in S$, there is a unique $g\in Z$ such that $gfx_2\in S$. By continuity, there is a neighborhood $Z$ such that every $x_1,x_2$ in the same atom of $\cC\vee\cD_Z$,
\[
\la\{f\in L|P(fx_1)\neq P(fx_2)\}<3\ep\la(L).
\]
\end{proof}

Since $rls(\cC\vee\cD|\cG)=rls(\cC|\cG)$, we can replace $\cC$ by $\cC\vee\cD$ and assume that each $C_i$ is a $(P,L,3\ep)$-ball. We show now that each $M_iC_j$ is a $(P,L,5\ep)$-ball. Indeed, if $x_1,x_2\in C_i$ and $g_1,g_2\in M_i$, then the distance between the $(L,P)$-names of $g_1x_1$ and $g_2x_2$ is
\[
\la\{f\in L | P(fg_1x_1)\neq P(fg_2x_2)\}\leq\la\{f\in L| P(fg_1x_1)\neq P(fg_2x_1)\} + 
\]
\[
+ \la\{f\in L|P(fg_2x_1)\neq P(fg_2x_2)\}.
\]
By the definition of $M_i$, we have that $fg_1=\de(f,g_1,g_2)fg_2$, where $\de(f,g_1,g_2)\in W$. Hence, the first summand is less than $\la\{f\in L|fg_2x_1\not\in B\}<\ep\la(L)$. Since $L$ is $(F,\ep)$-invariant, the second summand is less than $\ep\la(L)+\la\{f\in L| P(fx_1)\neq P(fx_2)\}<4\ep\la(L)$. We get that $h(G,X|\cG)\leq h(S|\cG_S)$.

For the other direction, we start with a definition.

\begin{defn} Given a partition $P$ of $S$ and $U\subset G$ open neighborhood of $1$, the $U$-fattening of $P$ is the partition $Q$ of $X$ into $|P|+1$ parts defined as follows: $Q(us)=P(s)$ if $u\in U, s\in S$, and $X\setminus US$ is the last part.
\end{defn}

Suppose we have a partition $\tP$ of $S$. Assume that $S$ is $U$-discrete for a neighborhood $U$ of the identity in $G$. Without loss of generality we can assume that $\la (U)<1$. Let $P$ be the $U$ fattening of $\tP$ to $X$. By definition, there is a F\o lner set $F$ and a partition $\cC =\{ C_1 ,\dots ,C_N \}$of $1-\ep \la (U)$ of $X$ such that for every $x\in X$, the set of $i$'s such that $E(1_{C_i}|\cG )(x)>0$ has size less than $\exp ((h+\ep )\la (F)$, and for every two points $x,y$ in the same part
\[
\la ( \{ f\in F | P(fx)\neq P(fy)\} )<\ep \la (F) .
\]
We can assume that $F$ is $(U,\ep )$-invariant and get that for every $u\in U$,
\[
\la ( \{ f\in F | P(fux)\neq P(fuy)\} )< 2\ep \la (F) .
\]
By Fubini, there is $u\in U$ such that $\cC$ covers more than $1-\ep$ of $uS$. It is now easy to check that $u^{-1}C_i \cap S$ is a collection of $(P,F,3\ep )$-balls. Since we did not increase its size, Theorem \ref{thm:SM} implies that $h(S|\cG _S)\leq h(G,X|\cG)+3\ep$.

\end{proof}

\section{CPE actions} \label{sec:cpe.actions}

We recall some definitions.

\begin{defn} Let $G$ be an amenable group with zero self entropy. An action of $G$ on a probability space $X$ is called completely positive entropy (or CPE for short) if, for any non-trivial partition $P$ of $X$, the entropy $h(G,P)$ is positive.
\end{defn}

\begin{defn} \label{unifmix} An action of an amenable group $G$ on a probability space $X$ is called uniformly mixing if, for every partition $P$ of $X$ and any $\ep >0$, there is a compact $K\subset G$ such that for any finite set $F\subset G$ that is $K$-separated (i.e. $gh^{-1}\not \in K$ for any two distinct $g,h\in F$), one has
\begin{equation} \label{eq:unif.mix}
\left | \frac{1}{|F|}H\left ( \bigvee \limits _{g\in F} gP \right ) -H(P)\right | < \ep .
\end{equation}
\end{defn}

\begin{defn} Let $G$ be an amenable group acting on a probability space $X$. The spectrum of this action is the associated $G$-representation on $L_2(X)$ given by
\[
(\rho(g)f)(x)=f(g^{-1}x).
\]
The spectrum is called Lebesgue with multiplicity $N$ (which can be infinity) if $L_2 (X)$ decomposes into direct sum of $N$ copies of the regular representation of $G$.
\end{defn}

\begin{thm} \label{CPEthm} Let $G$ be an amenable group with zero self entropy, and let $X$ be a probability space. Suppose that $G$ acts on $X$ in a CPE manner. Then
\begin{enumerate}
\item The action is uniformly mixing.
\item The spectrum of the action is Lebesgue with multiplicity $\aleph _0$.
\end{enumerate}
\end{thm}

\begin{proof} We first relativize the notions of CPE, uniform mixing, and spectrum as follows:

\begin{defn} Let $\cG$ be an invariant sub-sigma-algebra. The action of $G$ on $X$ is said to be relatively CPE over $\cG$ if, for any partition $P$, which is not measurable with respect to $\cG$, we have $h(G,P|\cG )>0$. We say that the action is uniformly mixing relative to $\cG$ if the inequality in \ref{unifmix} holds after we replace all entropies with relative entropies (with respect to $\cG$). Finally, the relative spectrum of the action is the $G$ representation $L_2(X,\cB )\ominus L_2(X,\cG )$ (by which we mean the orthogonal complement to $L_2 (X,\cG )$ in $L_2 (X,\cB )$).
\end{defn}

Note also that the notion of CPE makes sense also for cross-sections, since (relative) entropy is defined for them.

\begin{prop} \label{CPE:implies:CPE} With the same assumptions as in Theorem \ref{Abramov}, if the $G$ action on $X$ is completely positive entropy relative to $\cG$, then $S$ is completely positive entropy relative to $\cG$.
\end{prop}

\begin{proof} This follows from the proof of the second direction of Theorem \ref{Abramov}
\end{proof}

\begin{lem} There is a sequence of invariant sub-sigma-algebras $\dots \cG _2 \subset \cG _1 \subset \cG _0=\cB$ such that the action of $G$ on $(X,\cG _n)$ is CPE relative to $\cG _{n+1}$, and $\cap \cG _n$ is trivial (i.e. contains only null and conull sets).
\end{lem}

\begin{proof} We construct the sequence by induction. Suppose $\cG _n$ has been constructed. Choose a partition $P_n$ of $X$ that is measurable with respect to $\cG _n$, and such that $H(P)<\min\{h(G,\cG_n),\frac{1}{n}\}$. Let $\cF _n$ be the invariant sigma-algebra generated by $P$. The factor $(X,\cB )\to (X,\cF _n)$ need not be CPE, but it has a Pinsker factor, which is the sigma-algebra generated by all partitions $Q$ such that $h(Q | \cF _n)=0$. We take $\cG _{n+1}$ be this sigma-algebra. The first requirement on $\cG _n$ is satisfied, and the second follows since if $P$ is a partition measurable with respect to $\cap \cG _n$ then
\[
h(G,P) \leq h(G,\cG _{n}) = h(G,\cF _n)<\frac{1}{n} .
\]
So $h(G,P)=0$, and by the CPE assumption $P$ is trivial.
\end{proof}

For the second claim of the theorem, it is enough to show that $L_2 (\cG _n)\ominus L_2 (\cG _{n+1})$ is isomorphic to a countable sum of regular representations of $G$. We show that an analogue statement is also true for the first claim. First, recall the definition of the entropy of a partition, relative to a sigma-algebra.

\begin{defn} Let $P=\{ P_1,\ldots ,P_n \}$ be a partition of a probability space $(X,\cB ,m)$, and let $\cG \subset \cB$ be a sub-sigma-algebra. Let $p_i(x)=E(1_{P_i}|\cG )(x)$. Then for almost all $x$ we have that $p_i(x)\geq 0$ and $\sum p_i (x)=1$. We define
\[
H(P|\cG )=\int _X p_1(x)\log p_1 (x) + \ldots p_n(x)\log p_n(x) dx .
\]
\end{defn}

\begin{lem} Suppose that the action of $G$ on $X$ is uniformly mixing relative to all $\cG _n$'s. Then it is uniformly mixing.
\end{lem}

\begin{proof} Let $P$ be a partition and let $\ep >0$. Choose $n$ such that $H(P|\cG _n)>H(P)-\ep$. By the assumption, there is a compact $K\subset G$ such that if $F\subset G$ is $K$-separated then 
\[
\frac{1}{|F|}H\left ( \bigvee \limits _{g\in F}gP \right ) \geq \frac{1}{|F|}H \left ( \bigvee \limits _{g\in F} gP \Bigg |\cG _n \right ) > H(P| \cG _n )-\ep > H(P)-2\ep .
\]
But the inequality in the other direction is trivial.
\end{proof}

It is thus enough to show both relative claims for $n=1$. Choose a cross-section $S$ which is measurable with respect to $\cG$ and is $U$-discrete for some neighborhood of the identity in $G$. The cross-section $(S,\cB ,\mu ,\cR ,\al :\cR \to G)$ is such that $\al$ is $\cG$-measurable. By Proposition \ref{CPE:implies:CPE} this cross-section is CPE relative to $\cG$. By a theorem of \cite{CFW} there is a transformation $T:S\to S$ that is $\cG$-measurable and generates $\cR$. This gives us another cross-section $(S,\cB ,\mu ,\cR ,\be :\cR \to \bZ )$, and by Theorem \ref{transference} it is CPE relative to $\cG$. Again, by Proposition \ref{CPE:implies:CPE} we conclude that $T$ is CPE relative to $\cG$, hence uniformly mixing. We wish to transfer this back to $X$.

\begin{defn} We say that a cross-section $(S,\cB ,\mu ,\cR ,\al )$ is uniformly mixing over a sub-sigma-algebra $\cG$ (such that $\al$ is $\cG$-measurable) if, for all partitions $P$ of $S$ and $\ep >0$, there is compact $K\subset G$ such that for any sequence $\phi _i :S\to S$, $i=1,\dots ,n$ that satisfies

\begin{enumerate}

\item $(x,\phi _i(x))\in \cR$ for all $i$ and $x$.

\item $\al (x, \phi _i (x))\not \in K$.

\item $\phi _i$ is an isomorphism of measure spaces.

\end{enumerate}

we have that

\[
\left | \frac{1}{n} H \left ( \bigvee \limits _{i=1} ^{n} \phi _i P \Bigg | \cG \right )-H(P|\cG ) \right | <\ep .
\]
\end{defn}

In \cite{RW} it is shown that CPE for a $\bZ$ action implies uniform mixing for the cross-section it generates. Moreover, it shown there that (in our notation) the cross-section $(S,\cB ,\mu ,\cR ,\al )$ is uniformly mixing over $\cG$. It remains to show that this condition implies that the original $G$ action is uniformly mixing.

Suppose $P$ is a partition of $X$ and $\ep >0$. Choose a compact set $M\subset G$ such that $m(M\cdot S)>\ep$ and a neighborhood $V\subset G$ of the identity such that $S$ is $V$-discrete and for any $g\in V$,
\[
m(\{ x| P(x)\neq P(gx)\})<\ep .
\]
Finally, choose a finite set $E\subset M$ such that $M\subset EV$. Define a partition $Q$ on $S$ to be the common refinement of $gP$ where $g\in E$. Let $K\subset G$ be the compact set that one gets from the uniform relative mixing for the partition $Q$ and for $\ep \la (V)$. We claim that if $F\subset G$ is $KM$-separated then the inequality (\ref{eq:unif.mix}) in Definition \ref{unifmix} holds. Assuming the contrary, since the sets $gV\cdot S$, $g\in E$ cover almost all $X$, there is $g\in E$ such that the $fP$, $f\in F$ are not $\ep \la$ independent on $gV\cdot S$. We can also assume that $m(g_i gV\cdot S \cap MS)\geq (1-\ep )\la (V)$. Define $\phi _i:S\to S$ by $\phi _i (x)=hg_i gx$ if $g_i gx\in hS$ and $h\in M$ (and let $\phi _i (x)=x$ otherwise). Since the $g_iP$ are sufficiently dependent, we get that $\phi _i (Q)$ are also dependent. A contradiction.

We now move on to the spectral claim. We have that the action of $T$ on $S$ is CPE relative to $\cG$. By Rohlin-Sinai there is a sigma-algebra $\cH$ on $S$ such that $T\cH \subset \cH$, $\vee T^n \cH =\cB$, and $\wedge T^n \cH =\cG$. Choose an orthonormal basis $\phi _1 ,\phi _2 ,\dots$ to $L_2(\cH )\ominus L_2 (T\cH )$. For any $i$, every neighborhood $V\subset G$ of the identity such that $S$ is $V^2$-discrete, and every $g\in G$, define a function $\Phi _i ^{g,V}\in L_2 (X)$ by $\Phi _i ^{g,V}(x)=\phi _i(y)$ if $y\in S$ and $x\in gV\cdot y$, and $\Phi _i ^{g,V}(x)=0$ otherwise. Note that if such $y$ exists then it is unique.

\begin{lem} Assume that $S$ is $W$-discrete for some neighborhood of $1$ in $G$. Suppose that $V,U\subset G$ are such that $V^2,U^2 \subset W$. Then

\begin{enumerate}

\item If $i\neq j$ then $\Phi _i ^{g,V}$ and $\Phi _j ^{h,U}$ are orthogonal.

\item $<\Phi _i ^{g,V} , \Phi _i ^{h,U}>=\la (gV \cap hU)$.

\end{enumerate}

Moreover, the set $\{ \Phi _i ^{g,V} \}$ spans $L_2 (\cB )\ominus L_2 (\cG)$.
\end{lem}

\begin{proof} We want to compute the inner product of $\Phi _i ^{g,V}$ and $\Phi _j ^{h,U}$. Without loss of generality we may assume $h=1$. By definition, $\Phi _i ^{g,V} (x)\Phi _j ^{1,U}(x)=0$ unless $x\in gVS$ and $x\in US$. Hence,
\[
<\Phi _i ^{g,V},\Phi _j ^{h,U}>=\int _{s\in S} \int _{u\in U}\Phi _i ^{g,V} (us) \phi _j (s) 1_{\{ us\in gVS\} }du ds = 
\]

\[
\int _S \int _{U\cap gV} \Phi _i ^{g,V}(us) \phi _j(s) du ds+ \int _S \int _{\{ u\in U \setminus gV | us\in gVS\} }  \Phi _i ^{g,V}(us)\phi _j (s) du ds.
\]
In the first integral, $\Phi _i ^{g,V}(us)=\phi _i (s)$, so the integral is
\[
\int _S \phi _i (s)\phi _j (s) ds\int _{U\cap gV}du = <\phi _i , \phi _j>\la (U\cap gV).
\]
As for the second integral, for each $s\in S$ and $u\in U$, the value of $\Phi _i ^{g,V}(us)$ is either zero or equals $\phi _i (T^n)$ for some $n=n(s,u)\in \bZ \setminus \{ 0\}$. It is easy to see that $n$ is independent of $u$, and that the function $s\to n(s)$ is $\cG$-measurable. Moreover, the function $R$ that assigns to $s$ the measure of the set of $u\in U$ such that $us\in gVS$, is also $\cG$-measurable. The second integral equals
\[
\int _S \phi _i (s) \phi _j (T^{n(s)}s) R(s) ds = \sum _{k\neq 0} \int _{S}\phi _i(s) \phi _j (T^k s) R(s) 1_{\{ n(s)=k\} } ds.
\]
Suppose $k<0$. Then $s\mapsto \phi _j (T^k s)R(s) 1_{\{ n(s)=k \}}$ is $T^{-1}\cH$-measurable, and since $\phi _i \in L_2(\cH )\ominus L_2(T^{-1}\cH )$, the integral is zero. If $k>0$, the same argument works if we replace $i$ and $j$. This proves 1. and 2. above.

Finally, assume $\Psi \in L_2 (X)$. Then for almost all $g\in G$ we have that $\psi _g = g\Phi |_{S}$ is in $L_2(S)$. Suppose $\Psi \perp \Phi _i ^{g,V}$ for all $i$, $g$ and $V$. Then $\psi _g \perp T^n\phi _i$ for all $i$ and $n$. Hence for almost all $g\in G$, we have $\psi _g \in L_2(S,\cG )$. Since $\cG$ is $G$-invariant, $\Psi \in L_2(X,\cG )$.
\end{proof}

By the lemma, for every $i$, the map $1_{gV} \mapsto \Phi _i ^{g,V}$ extends to a $G$-equivariant  isometric embedding $L_2 (G) \to L_2 (X,\cB )\ominus L_2 (X,\cG )$, and the images of those maps for different $i$'s are orthogonal. By the second assertion of the lemma, these copies span.

\end{proof}

\begin{footnotesize}
\begin{quote}

Nir Avni\\
Department of Mathematics\\
Harvard University\\
Cambridge, MA 02138\\
{\tt avni.nir@gmail.com}

\end{quote}
\end{footnotesize}

\end{document}